\documentclass[11pt]{article}
\usepackage{amsfonts,latexsym}

\title{Non-abelian Reciprocity Laws on a Riemann Surface}
\author{Ivan Emilov Horozov}
\date{February 16, 2010}

\newcommand{\beq}{\begin{equation}}
\newcommand{\eeq}{\end{equation}}
\newcommand{\beqa}{\begin{eqnarray}}
\newcommand{\eeqa}{\end{eqnarray}}
\newcommand{\beaa}{\begin{eqnarray*}}
\newcommand{\ben}{\begin{eqnarray*}}
\newcommand{\eaa}{\end{eqnarray*}}
\newcommand{\een}{\end{eqnarray*}}

\newcommand{\text}{\textrm}
\newcommand \nc {\newcommand}
\nc \proof {\noindent {\em{Proof.\/ }}} \nc \qed {$\Box$\hfill}
\newtheorem{theorem}{Theorem}[section]
\newtheorem{lemma}[theorem]{Lemma}
\newtheorem{proposition}[theorem]{Proposition}
\newtheorem{corollary}[theorem]{Corollary}
\newtheorem{definition}[theorem]{Definition}
\newtheorem{example}[theorem]{Example}
\newtheorem{remark}[theorem]{Remark}
\newtheorem{conjecture}[theorem]{Conjecture}
\newtheorem{question}[theorem]{Question}
\nc \bth[1] {\begin{theorem}\label{t#1} } \nc \ble[1]
{\begin{lemma}\label{l#1} } \nc \bpr[1]
{\begin{proposition}\label{p#1} } \nc \bco[1]
{\begin{corollary}\label{c#1} } \nc \bde[1]
{\begin{definition}\label{d#1}\rm } \nc \bex[1]
{\begin{example}\label{e#1}\rm } \nc \bre[1]
{\begin{remark}\label{r#1}\rm } \nc \bcon[1]
{\begin{conjecture}\label{con#1}\rm } \nc \bque[1]
{\begin{question}\label{que#1}\rm }
\nc {\eth} { \end{theorem} } \nc {\ele} { \end{lemma} } \nc
{\epr}{ \end{proposition} } \nc {\eco} { \end{corollary} } \nc
{\ede} {\end{definition} } \nc {\eex} { \end{example} } \nc {\ere}
{\end{remark} } \nc {\econ} { \end{conjecture} } \nc {\eque}
{\end{question} }
\nc \eqref[1] {{\rm{(\ref{#1})}}} \nc \thref[1]{Theorem \ref{t#1}}
\nc \leref[1]{Lemma \ref{l#1}} \nc \prref[1]{Proposition
\ref{p#1}} \nc \coref[1]{Corollary \ref{c#1}} \nc
\deref[1]{Definition \ref{d#1}} \nc \exref[1]{Example \ref{e#1}}
\nc \reref[1]{Remark \ref{r#1}} \nc
\conref[1]{Conjecture\ref{con#1}}

\def\a{\alpha}
\def\b{\beta}

\nc \Wr {Wr} \nc \GRN { \Gr^{(N)} }
\nc \GRA[1] { \Gr_A^{(#1)} }   
\nc \GRAN { \GRA{N} } \nc \GrA[1] { \Gr_A(#1) }\nc \GrAa {
\GrA{\alpha} }
\nc \GRB[1] { \Gr_B^{(#1)} }   
\nc \GRBN { \GRB{N} } \nc \GrB[1] { \Gr_B(#1) } \nc \GrBb {
\GrB{\beta} }
\nc \GRMB[1] { \Gr_{MB}^{(#1)} }   
\nc \GRMBN { \GRMB{N} } \nc \GrMB[1] { \Gr_{MB}(#1) } \nc \GrMBb {
\GrMB{\beta} }

\begin{document}

\title{{\LARGE\bf{Non-abelian Reciprocity Laws on a Riemann Surface}}}

\author{
I. ~Horozov
\thanks{E-mail: horozov@uni-tuebingen.de}
\\ \hfill\\ \normalsize \textit{Mathematisches Institut,}\\
\normalsize \textit{Universit\"at T\"ubingen, Auf der Morgenstelle 10,}\\
\normalsize \textit {72076 T\"ubingen, Germany }  \\
}
\date{}
\maketitle

\begin{abstract}
On a Riemann surface there are relations among the periods of
holomorphic differential forms, called Riemann's relations. If one
looks carefully in Riemann's proof, one notices that he uses
iterated integrals. What I have done in this paper is to
generalize these relations to relations among generating series of
iterated integrals. Since the main result is formulated in terms
of generating series, it gives infinitely many relations - one for
each coefficient of the generating series. The lower order terms
give the well known classical relations.  The new result is
reciprocity for the higher degree terms, which give non-trivial
relations among iterated integrals on a Riemann surface. As an
application we refine the definition of Manin's non-commutative
modular symbol in order to include Eisenstein series. Finally, we
have to point out that this paper contains some constructions
needed for multidimensional reciprocity laws like a refinement of
one of the Kato-Parshin reciprocity laws.
\end{abstract}


\tableofcontents
\setcounter{section}{-1}
\section{Introduction}
This paper is the first one from a series of papers on reciprocity
laws on a complex variety based on properties of a fundamental
group. Here we limit ourselves to the case of Riemann surfaces. We
obtain infinitely many reciprocity laws on a Riemann surface,
which we write as a reciprocity law for a generating series. The
first few reciprocity laws are known, but the rest are new. The
second paper from this series gives a reciprocity law for a new
symbol, which is a refinement of Parshin symbol (see \cite{P2},
\cite{P3}) for a complex surface.

In both papers we use properties of a certain fundamental group by
examining iterated integrals in the sense of K.-T. Chen (see
\cite{Ch}).

We define a generating series of iterated integrals $F$ with good
analytic properties. And we prove a reciprocity law for such a
generating series $F$. Since we have a generating series, each
coefficient gives us a reciprocity law. The first few terms give
us known reciprocity laws. For a term of degree $3$, we obtain a
new reciprocity law \cite{H}, which we use in the next paper on
refinement of one of the Kato-Parshin (\cite{Ka}) reciprocity laws
for the Parshin symbol (\cite{P2}, \cite{P3}). The lower degree
terms give well known reciprocity laws. In degree $1$ we recover
the following reciprocity: the sum of the residues of a
differential form on a Riemann surface is zero. In degree $2$, we
recover Riemann relations among differential forms of the $3$-rd
kind \cite{GH}. Also from degree $2$, we recover Weil reciprocity.
If one looks carefully at a proof of Riemann's relations (see
\cite{GH}), one notices that they use iterated integrals. What I
have done in this paper is to generalize these relations to
relations among generating series of iterated integrals $F$. Since
it is formulated in terms of generating series, it gives
infinitely many relations - one for each coefficient of the
generating series. This paper establishes the main tools, which
will be used in papers that establish a new reciprocity law for
portion of Kato-Parshin symbols. For a complex surface see
\cite{H}. For refinement of Parshin symbol for higher dimensional
complex varieties, we have the a prove for the new reciprocity
law. This in the $3$-rd paper from the series, which is in
preparation.

We construct the generating series of iterated integrals $F$ as a
solution of an ordinary differential equation, following an idea
of Manin \cite{M}. For each closed loop $\sigma$ on a Riemann
surface, we consider the generating series of the iterated
integrals $F_{\sigma}$. On a Riemann surface there is "good"
choice of generators of the fundamental group, with only one
relation. For each of the generators of the fundamental group we
consider the generating series of iterated integrals $F_{\sigma}$.
Composition of paths corresponds to a composition of generating
series of iterated integrals. For two loops $\sigma_1$ and
$\sigma_2$ with common base point, we have
$F_{\sigma_1}F_{\sigma_2}=F_{\sigma_1\sigma_2}$ for the generating
series. The relation among the generators in the fundamental group
gives us a relation among the generating series of iterated
integrals $F_{\sigma}$ for various loops $\sigma$. We call this
relation a reciprocity law. This is the key idea in the
construction of the reciprocity law.

Before we formulate the final statement of the reciprocity law, we
consider several simplifications. The first one is for a simple
loop around a pole of any of the differential forms. We present
constructively the process of taking residues of iterating
integrals, which we generalize to taking residues of generating
series of iterated integrals. Another simplification that we make
is for $F_{[\a_i,\b_i]}$, where $[\a_i,\b_i]$ is commutator of
$\a_i$ and $\b_i$ loops on a Riemann surface. When all this is
done, we formulate the reciprocity law.

In the last section we apply the developed ideas about
non-commutative reciprocity law to Manin's non-commutative modular
symbol. We are able to extend his construction  so that we can
include Eisenstein series.

Where is the non-abelian group in the reciprocity law? An iterated
integral over a loop depends on the loop. Since iterated integrals
are homotopy invariant, we have that an iterated integral over a
loop depends only on the element of the fundamental group that
this loop represents. How much do iterated integrals distinguish
elements of the fundamental group? I thank one of the referees,
who pointed out that conjecture by Parshin \cite{P1} that iterated
integrals of differential forms of the third kind on a Riemann
surface capture precisely the pro-unipotent part of the
fundamental group. For some progress and intuition in this
direction one might look at \cite{Q}. In terms of nilpotent
variations of mixed Hodge structures related to the fundamental
group, one can look at the following papers \cite{G} and
\cite{DG}.

There are non-commutative reciprocity laws by Brylinski and
McLaughlin (see \cite{BrMc}) on a complex variety. There
reciprocity laws are non-abelian, since the authors consider not
only a classifying space $BGL$ of the structure sheaf, but also a
certain central extension. In this paper the reciprocity laws are
non-abelian for a different reason. The reciprocity laws are based
on properties of the fundamental group on a punctured Riemann
surface $X$, which is a non-abelian group in most of the
interesting situations. Also, with the use of iterated integrals,
we capture non-abelian quotient of the fundamental group
$\pi_1(X,P)$, which is strictly larger than the maximal abelian
quotient $H_1(X)$.

I would like to make one remark about the reciprocity laws in this
paper. Instead of working on a Riemann surface, one can consider a
smooth algebraic curve over the algebraic closure of the rational
numbers. Using the construction in \cite{G} and \cite{DG} we can
conclude that the iterated integrals of algebraic differential
forms with logarithmic poles over a path give periods in the sense
of algebraic geometry. So the reciprocity laws that we describe in
this paper are relations among periods.

We are going to use explicit integrals. However, in a more
theoretical approach, one can consider framed mixed Hodge
structures associated to the integrals that we consider. The
process of taking residues of iterated integrals corresponds to
taking co-product of the corresponding framed mixed Hodge
structure (see \cite{G}).

The non-abelian reciprocity law in this paper has a generalization
to higher dimensions. Even in dimension 2, new phenomenons occur:
instead of reciprocity law for the Parshin symbol, we obtain a new
reciprocity law for a refinement of the Parshin symbol (see
\cite{H}).

\section{Background on iterated integrals}
This section establishes both the notation and the main properties
of iterated integrals, which we are going to use throughout the
paper. We recall well known properties of iterated integrals,
which we are going to use heavily in the rest of the paper. One
can look at K.-T. Chen \cite{Ch} and Goncharov \cite{G} for proofs
of the statements and for more properties of iterated integrals.
For the differential equation (\ref{DiffEq}), we follow the idea
of Manin \cite{M}.
\subsection{Definition of iterated integrals}
\begin{definition}
\label{def1} Let $\omega_1,\dots,\omega_n$ be holomorphic 1-forms
on a simply connected open subset $U$ of the complex plane
${\mathbb C}$. Let
$$\gamma:[0,1]\rightarrow U$$ be a path.
We define an iterated integral of the forms
$\omega_1,\dots,\omega_n$ over the path $\gamma$ to be
$$\int_\gamma \omega_1\circ\dots\circ\omega_n=
\int\dots\int_{0\leq t_1\leq\dots\leq t_n\leq 1}
\gamma^*\omega_1(t_1)\wedge\dots\wedge\gamma^*\omega_n(t_n).$$
\end{definition}

It is called iterated because it can be defined inductively by
$$\int_\gamma \omega_1\circ\dots\circ\omega_n=
\int_0^1
(\int_{\gamma|[0,t]}\omega_1\circ\dots\circ\omega_{n-1})\gamma^*\omega_n(t).$$

\subsection{Homotopy invariance of iterated integrals}

\begin{theorem} \label{thm1} Let $\omega_1,\dots,\omega_n$ be holomorphic 1-forms on
a simply connected open subset $U$ of the complex plane ${\mathbb
C}$. Let
$$H:[0,1]\times [0,1]\rightarrow U$$ be a homotopy, fixing the end
points, of paths
$$\gamma_s:[0,1]\rightarrow U$$ such that $\gamma_s(t)=H(s,t)$,
and for fixed $s$, we have a path $\gamma_s:[0,1]\rightarrow U$.
Then
$$\int_{\gamma_s} \omega_1\circ\dots\circ\omega_n$$
is independent of $s$. 
\end{theorem} 

\subsection{Differential equation}

When we consider an iterated integral, we can let the end point
vary in a small neighborhood. Then the iterated integral becomes
an analytic function.

Let $\omega_1,\dots,\omega_n$ be differentials of the 3rd kind on
a Riemann surface $X$. Following an idea of Manin \cite{M}, we
consider the differential equation
\begin{equation} dF=F\sum_{i=1}^n A_i\omega_i,\label{DiffEq}
\end{equation}
where $A_1,\dots,A_n$  are non-commuting formal variables. Let $P$
be a point on $X$ such that none of the differential forms has a
pole at $P$. It is easy to check that the solution of $F(z)$ with
initial condition $F(P)=1$ is
$$F(z)=1+\sum_{i}(A_i\int_P^z \omega_i) +
\sum_{i,j}A_iA_j\int_P^z\omega_i\circ\omega_j +
\sum_{i,j,k}A_iA_jA_k\int_P^z\omega_i\circ\omega_j\circ\omega_k+\dots.$$
The summation is over all iterated integrals of the given
differential forms. Note that
$$d\int_P^z \omega_i\circ\dots\circ\omega_j\circ\omega_k =
(\int_P^z \omega_i\circ\dots\circ\omega_j)\omega_k.$$

\subsection{Multiplication formulas}
We can take a path $\gamma$ from $P$ to $z$. We denote the
solution of the differential equation by $F_{\gamma}$. If
$\gamma_1$ is a path that ends at $Q$ and $\gamma_2$ is a path
that starts at $Q$ we can compose them. Denote the composition by
$\gamma_1\gamma_2$. 
\begin{theorem} \label{thm2} (Composition of paths) With the above notation, we have
$$F_{\gamma_1}F_{\gamma_2}=F_{\gamma_1\gamma_2}.$$
\end{theorem}
\begin{corollary} \label{co1}(Composition of paths) Let $\omega_1,\dots,\omega_n$ be
differential forms, some of them could repeat. Let also $\gamma_1$
be a path that ends at $Q$ and $\gamma_2$ be a path that starts at
$Q$. We can compose them. Denote the composition by
$\gamma_1\gamma_2$. Then
$$\int_{\gamma_1\gamma_2} \omega_1\circ\dots\circ\omega_n=
\sum_{i=0}^n \int_{\gamma_1}\omega_1\circ\dots\circ\omega_i
\int_{\gamma_2} \omega_{i+1}\circ\dots\circ\omega_n,$$ where for
$i=0$ we define
$\int_{\gamma_1}\omega_1\circ\dots\circ\omega_i=1$, and similarly,
for $i=n$ we define $\int_{\gamma_2}
\omega_{i+1}\circ\dots\circ\omega_n=1$.
 \end{corollary}
\subsection{Shuffle relations}
\begin{definition}
\label{def2} Denote by $Sh(m,n)$ the shuffles, which are
permutations $\tau$ of the set $\{1,\dots,m,m+1,\dots,m+n\}$ such
that
$$\tau(1)<\tau(2)<\dots<\tau(m)$$ and
$$\tau(m+1)<\tau(m+2)<\dots<\tau(m+n).$$
\end{definition}

\begin{theorem} \label{thm3}(Shuffle relation) Let
$\omega_1,\dots,\omega_m,\omega_{m+1},\dots,\omega_{m+n}$ be
differential $1$-forms, some of them could repeat. Let also
$\gamma$ be a path that does not pass through any of the poles of
the given differential forms. Then
$$\int_{\gamma} \omega_1\circ\dots\circ\omega_m\int_{\gamma}
\omega_{m+1}\circ\dots\circ\omega_{m+n}=
\sum_{\tau\in Sh(m,n)}\int_{\gamma}
\omega_{\tau(1)}\circ\omega_{\tau(2)}\circ\dots\circ\omega_{\tau(m+n)}.$$
 \end{theorem}
\subsection{Reversing the path}

\begin{lemma} \label{lem1}(Reversing the path) Let $\gamma$ be a path. Let
$\gamma^{-1}$ be the same path but going to the opposite
direction. Then
$$\int_{\gamma^{-1}}\omega_1\circ\omega_2\circ\dots\circ\omega_n=
(-1)^n\int_{\gamma}\omega_n\circ\omega_{n-1}\circ\dots\circ\omega_1.$$
 \end{lemma}
\section{Non-abelian reciprocity law on a Riemann surface}
This section is the heart of the article. It continuous with the
topic of iterated integrals. However, now this is done in the
direction of building the reciprocity law. It is more technical
that the previous section.And it ends with the statement of the
non-abelian reciprocity law on Riemann surfaces.
\subsection{Iterated integrals over a loop around a pole}
Let $U$ be an open simply connected subset of the complex plane.
We can assume that $0$ belongs to $U$. Let $f_1dz,\dots,f_ndz$ be
holomorphic differentials on $U$. We are going to iterate these
differential forms together with the form $dz/z$. Let $P$ be a
point in $U$ different from $0$. And let $\sigma$ be a simple loop
in $U$ that starts and ends at $P$ and goes around $0$ once in a
counterclockwise direction. We can assume that $\sigma$ does not
intersect itself. Let $\gamma_{\epsilon}$ be a path that starts at
$P$ and ends at $\epsilon$ for a point $\epsilon\neq 0$ in $U$
very close to $0$. For convenience we take $\epsilon$ to be a
positive real number. Let $\gamma$ be a path starting from $P$ and
ending at $0$, which is the limit of $\gamma_{\epsilon}$, when
$\epsilon$ tends to zero. We define also $\sigma_\epsilon$ to be a
loop that starts and ends at $\epsilon$ and goes around $0$ once
in a counterclockwise direction along a circle of radius
$\epsilon$. We can deform $\sigma$ homotopic in $U-\{0\}$ to
$\gamma_\epsilon\sigma_\epsilon\gamma_\epsilon^{-1}$. Note that
iterated integrals are invariant with respect to homotopic
deformation of the path of integration (see theorem \ref{thm1}).
We are going to use corollary \ref{co1} for the composition of the
paths $\gamma_\epsilon\sigma_\epsilon\gamma_\epsilon^{-1}$.

\begin{lemma} \label{lem2}With the above notation
$$\int_{\sigma_\epsilon}\frac{dz}{z}\circ\dots\circ\frac{dz}{z}=\frac{(2\pi
i)^r}{r!},$$ where we iterate the form $dz/dz$ with itself
$r$-times \end{lemma}

\proof We can take the following parametrization of $z$ along
$\sigma_\epsilon$. Let $z=\epsilon e^{2\pi it}$ for $0\leq t \leq
1$. Then $dz/z=2\pi idt$. Therefore,
$$\int_{\sigma_\epsilon}\frac{dz}{z}\circ\dots\circ\frac{dz}{z}=
(2\pi i)^r\int_0^1 dt\circ\dots\circ dt=\frac{(2\pi i)^r}{r!}.$$

\begin{corollary} \label{co2} If $\omega_1,\dots,\omega_r$ are holomorphic forms on
$U-\{0\}$ with simple poles at $0$ then
$$\int_{\sigma_\epsilon}\omega_1\circ\dots\circ\omega_r=
\frac{1}{r!}\prod_{i=1}^r Res_0 \omega_i.$$ \end{corollary}

\begin{lemma} \label{lem3} Consider an iterated integral over $\sigma_\epsilon$ of the
differential forms $$\frac{dz}{z}, f_1dz, \dots, f_ndz$$ in
various order with possible repetition, so that not all of them
need to appear in the integral. If at least one of the holomorphic
differentials $f_1dz,\dots,f_ndz$ on $U$ is present in the
integral then the limit of the iterated integral, as $\epsilon$
goes to
$0$, is $0$. \end{lemma} 

\proof Expand the holomorphic differential forms around $z=0$.
Consider the following parametrization of the variable $z$ along
the path $\sigma_\epsilon$: $z=\epsilon e^{2\pi it}$ for $0\leq t
\leq 1$. It is enough to prove the statement of the lemma for
$f_i(z)=z^{n_i}$. (The general statement will follow since we have
an uniform convergence of the power series of $f_i$ on compact
subsets of $U$.) Note that $dz/z=2\pi idt$ and $z^ndz=2\pi i
\epsilon^{n+1} e^{2\pi i(n+1)t}dt$ for some $n\geq 0$. Consider
the parametrization in terms of $t$ for $0\leq t\leq s$, where $s$
is close to $1$. Then the iterated integral over $0\leq t\leq s$
becomes
$$C\epsilon^N e^{2\pi i Ns}g(s),$$
where $C$ is a constant $N>0$ and $g(t)$ is a polynomial in $t$.
So the limit as $s$ approaches $1$ will be $C\epsilon^N g(1)$. So
the iterated integral has value  $C\epsilon^N g(1)$. Finally, this
value approaches zero as $\epsilon$ tends to zero, which proves
the lemma.

\begin{lemma} \label{lem4} Let $f_1dz,\dots,f_ndz$ be holomorphic forms on
$U-\{0\}$. Let $f_i$ and $f_{i+1}$ be also holomorphic at $0$.
Denote by $\omega^{\circ r}$ $r$-fold iteration of $\omega$ Then
$$\lim_{\epsilon\rightarrow 0}
\sum_{j=0}^r\int_{\gamma_\epsilon} f_1dz\circ\dots\circ
f_idz\circ\left(\frac{dz}{z}\right)^{\circ j}
\int_{\gamma_\epsilon^{-1}} \left(\frac{dz}{z}\right)^{\circ
(r-j)}\circ f_{i+1}dz\circ\dots\circ f_ndz=0.$$ The limit is still
zero in the cases when the set $\{f_1dz,\dots,f_idz\}$ is empty
and/or the set $\{f_{i+1}dz,\dots,f_ndz\}$ is empty. \end{lemma}

\proof Consider the integral $$\int_{\gamma_\epsilon}
f_1dz\circ\dots\circ f_idz$$ as a function of $\epsilon$. It is a
sum of terms of the type a constant times $\epsilon^k
log^l(\epsilon)$. Since $f_i$ is holomorphic, we have that if
$l\geq 0$ then $k\geq 0$. Also if $l>0$ then $k>0$, because the
integral is convergent for $\epsilon=0$. When $\epsilon$ tends to
zero $\epsilon^k log^l(\epsilon)$ tend to zero. Now, consider the
integral
$$\int_{\gamma_\epsilon} f_1dz\circ\dots\circ
f_idz\circ\left(\frac{dz}{z}\right)^{\circ j}.$$ As a function of
$\epsilon$ it is a sum of terms of the type constant times
$\epsilon^k log^{l+j}(\epsilon)$ for $k\geq 0$ and $l\geq 0$. From
the previous considerations if $l>0$ then $k>0$ and the term tends
to zero as $\epsilon$ approaches zero. Thus, the only term that
does not tend to zero as $\epsilon$ approaches zero is the
constant term of

$$\int_{\gamma_\epsilon}f_1dz\circ\dots\circ f_idz$$
times $log^i(\epsilon)$. Then
$$\lim_{\epsilon\rightarrow 0}\int_{\gamma_\epsilon} f_1dz\circ\dots\circ
f_idz\circ\left(\frac{dz}{z}\right)^{\circ j}
-\frac{1}{j!}log^j(\epsilon)\int_{\gamma} f_1dz\circ\dots\circ
f_idz=0.$$

For the other integral in this lemma we use
$$\int_{\gamma_\epsilon^{-1}} \left(\frac{dz}{z}\right)^{\circ
(r-j)}\circ f_{i+1}dz\circ\dots\circ
f_ndz=(-1)^{n-i+r-j}\int_{\gamma_\epsilon} f_ndz \circ\dots\circ
f_{i+1}dz\circ \left(\frac{dz}{z}\right)^{\circ (r-j)}.$$ Using
the same arguments as in the beginning of the proof we find that

$$\lim_{\epsilon\rightarrow 0}\int_{\gamma_\epsilon} f_ndz \circ\dots\circ
f_{i+1}dz\circ \left(\frac{dz}{z}\right)^{\circ (r-j)}-
\frac{1}{(r-j)!}log^{r-j}(\epsilon)\int_{\gamma} f_ndz
\circ\dots\circ f_{i+1}dz=0.$$

Note also that $$(-1)^{n-i+r-j}\int_{\gamma} f_ndz
\circ\dots\circ f_{i+1}dz=(-1)^{r-j}\int_{\gamma^{-1}} f_{i+1}dz
\circ\dots\circ f_{n}dz.$$

Therefore,
$$
\begin{tabular}{lll}
&$\lim_{\epsilon\rightarrow 0} \sum_{j=0}^r
\int_{\gamma_\epsilon}f_1dz\circ\dots\circ f_idz\circ
\left(\frac{dz}{z}\right)^{\circ j}\int_{\gamma_\epsilon^{-1}}
\left(\frac{dz}{z}\right)^{\circ
(r-j)}\circ f_{i+1}dz\circ\dots\circ f_ndz=$\\
\\
=&$\int_{\gamma} f_1dz\circ\dots\circ f_idz\int_{\gamma^{-1}}
f_{i+1}dz \circ\dots\circ
f_{n}dz\sum_{j=0}^r\lim_{\epsilon\rightarrow
0}\frac{1}{j!}(-1)^{r-j}\frac{1}{(r-j)!}log^r(\epsilon).$\\
\end{tabular}
$$
Finally, we have
$$\sum_{j=0}^r
\frac{1}{j!}log^j(\epsilon)(-1)^{r-j}\frac{1}{(r-j)!}log^{r-j}(\epsilon)=0,$$
using binomial coefficients, after multiplying by $r!$. Thus, the
limit in the lemma is zero.

Now we are ready to give the general algorithm for expressing an
iterated integral over a loop $\sigma$ around $0$ in terms of
integrals over a path $\gamma$ starting at $P$ and ending at $0$
and residues at $z=0$.

\begin{theorem} \label{thm4} Let $\omega_1,\dots,\omega_n$ be holomorphic forms on
$U$. Let $i_1,\dots,i_m$ be integers such that $$0\leq
i_1<i_2<\dots<i_m\leq n.$$ Let $j_1,\dots,j_m$ be positive
integers. Then
$$
\begin{tabular}{ll}
&$\int_{\sigma}\omega_1\circ\dots\omega_{i_1}\circ
\left(\frac{dz}{z}\right)^{\circ j_1}
\circ\omega_{i_1+1}\circ\dots\circ\omega_{i_2}\circ
\left(\frac{dz}{z}\right)^{\circ j_2}\circ\omega_{i_2+1}\circ\dots$\\
\\
&$\dots\circ\omega_{i_m}\circ\left(\frac{dz}{z}\right)^{\circ j_m}
\circ\omega_{i_m+1}\circ\dots\circ\omega_n=$\\
\\
=&$\sum_{k=1}^m \frac{(2\pi i)^{j_k}}{j_k!}
\int_\gamma\omega_1\circ\dots\circ\omega_{i_1}\circ
\dots\circ\left(\frac{dz}{z}\right)^{\circ j_{k-1}}
\circ\omega_{i_{k-1}+1}\circ\dots\circ\omega_{i_k}\times$\\
\\
&$\times\int_{\gamma^{-1}}\omega_{i_k+1}\circ\dots\circ\omega_{i_{k+1}}\circ
\left(\frac{dz}{z}\right)^{\circ
j_{k+1}}\circ\omega_{i_{k+1}+1}\circ
\dots\circ\omega_{i_m}\circ\left(\frac{dz}{z}\right)^{\circ j_m}
\circ\omega_{i_m+1}\dots\circ\omega_n.$
\end{tabular}
$$
\end{theorem} 
\proof First we use that $\sigma$ is homotopic to
$\gamma_\epsilon\sigma_\epsilon\gamma_\epsilon^{-1}$. Then we use
the formula for composition of paths in corollary \ref{co1},
expressing the integral over $\sigma$ in terms of a sum of
products of an integral over $\gamma_\epsilon$, an integral over
$\sigma_\epsilon$ and an integral over $\gamma_\epsilon^{-1}$.
Consider the portion of this sum, where there are no differential
forms integrated over $\sigma_\epsilon$. The sum of all such terms
is zero because by corollary \ref{co1} this is the same as the
integral over $\gamma_\epsilon\gamma_\epsilon^{-1}$, which is
homotopic to the constant loop at $P$.

Next, we examine what possible iterated integrals over
$\sigma_\epsilon$ we can have. If we have a holomorphic
differential form at zero in the iterated integral over
$\sigma_\epsilon$ then, by lemma \ref{lem3}, the value of the
integral tends to zero as $\epsilon$ approaches zero. Therefore,
it is enough to consider only iterations of $dz/z$ over
$\sigma_\epsilon$.

Consider iterated integrals of $dz/z$ over $\sigma_\epsilon$. If
the corresponding iterated integral over $\gamma$ ends with $dz/z$
and/or the corresponding iterated integral over
$\gamma_\epsilon^{-1}$ starts with $dz/z$ then by lemma \ref{lem3}
the sum of all such integrals tends to zero as $\epsilon$
approaches zero.

Therefore, the only terms we have to consider are the ones in the
theorem.

\subsection{Generating series of iterated integrals over a loop around a pole}
Consider $n$ differential 1-forms of with simple poles
$\omega_1,\dots,\omega_n$, defined on an open and simply connected
set $U\subset {\mathbb C}$. Consider the differential equation
(\ref{DiffEq})
$$dF=F\sum_{i=1}^n A_i\omega_i$$
where $A_1,\dots,A_n$ are non-commuting formal variables. Let
$Q\in U$ be a point where at least one of the differential forms
has a pole. Let $P\in U$ be a point, which is not a pole for any
of the differential forms. Consider a simple loop $\sigma$ with
the following properties: it starts at $P$; it does not
self-intersect; it bounds a region $V$ homeomorphic to a disk; and
the only pole of the differential forms that lies inside $V$ is
$Q$.

Let us simplify the solution $F_{\sigma}$, using theorem
\ref{thm4}. Define $\gamma$ to be the path, starting at $P$ and
ending at $Q$, that sits in the region $V$, bounded by $\sigma$.
Using theorem \ref{thm4}, we can decompose $F_{\sigma}$. Let
$F^{reg}_{\gamma}$ be regularization of $F_{\gamma}$, which
contains only the iterated integrals over the path $\gamma$, whose
iteration does not end with a differential form that has a pole at
$P$. The series $F_{\gamma}$ contains the summand $1$ and also all
the iterated integrals mentioned above times the corresponding
non-commuting variables. Similarly, $F^{reg}_{\gamma^{-1}}$ is
regularization of $F_{\gamma^{-1}}$, which contains only the
iterated integrals over the path $\gamma$, whose iteration does
not start with a differential form that has a pole at $Q$, where
$\gamma^{-1}$ is the reversed path
of $\gamma$. 
\begin{lemma} \label{lem5} With the notation in this subsection, we have
$$(F^{reg}_{\gamma})^{-1}=F^{reg}_{\gamma^{-1}}.$$
\end{lemma} \proof It follows from lemma \ref{lem1}. 

We call $F^{Res}_{\sigma}$ the residual part of $F_{\sigma}$,
defined by the portion of $F_{\sigma}$ that contains only the
iterated integrals of differential forms that have a pole at $Q$.
In particular, $F^{Res}_{\sigma}$ does contain the constant $1$.
Let
$$F^{Res+}_{\sigma}=-1+F^{Res}_{\sigma}$$
be the residual part without the constant term $1$. Recall, $Q$ is
a point, where at least one of the differential forms has a pole.
Then we have the following version of theorem \ref{thm4} in terms
of
generating series. 
\begin{theorem} \label{thm5} With the notation from this subsection, we have
$$F_{\sigma}=1+F^{reg}_{\gamma}F^{Res+}_{\sigma}F^{reg}_{\gamma^{-1}}.$$
\end{theorem} \proof From theorem \ref{thm4} we have
$$F_{\sigma}=F^{reg}_{\gamma}F^{Res}_{\sigma}F^{reg}_{\gamma^{-1}}.$$
We simplify the right hand side, using lemma \ref{lem5}
$$F^{reg}_{\gamma}F^{reg}_{\gamma^{-1}}=1.$$
$$F^{reg}_{\gamma}F^{Res}_{\sigma}F^{reg}_{\gamma^{-1}}=F_{\sigma}=
F^{reg}_{\gamma}(1+F^{Res+}_{\sigma})F^{reg}_{\gamma^{-1}}
=1+F^{reg}_{\gamma}F^{Res+}_{\sigma}F^{reg}_{\gamma^{-1}}.$$
\subsection{Generating series over $\a$ and $\b$ cycles of a Riemann surface}
Let $X$ be a Riemann surface of genus $g$. Let
$\omega_1,\dots,\omega_n$ be differential forms of the third kind
on $X$. Let also $\a_1,\b_1,\dots,\a_g,\b_g$ be loops on $X$
starting at $P$, which do not pass through a pole of any of the
differential forms, such that they generate the fundamental group
$\pi_1(X,P)$ with only one relation
$$[\a_1,\b_1]\dots[\a_g,\b_g]=1.$$
We used the notation $[\a,\b]=\a\b\a^{-1}\b^{-1}$.

We are going to simplify
$$F_{[\a,\b]}=F_{\a}F_{\b}F_{\a^{-1}}F_{\b^{-1}}.$$ Note that if we
take the constant term $1$ from $F_{\b}$ then we will have
$$F_{\a}F_{\a^{-1}}F_{\b^{-1}}=F_{\b^{-1}},$$
which follows from the homotopy invariance of $F_{\a\a^{-1}}$
(theorem \ref{thm1}. In order to capture such a cancellation, we
define
for a path $\gamma$ $$F^+_{\gamma}=F_{\gamma}-1.$$ 
\begin{lemma} \label{lem6} With the above notation
$$F_{[\a,\b]}=1+F^+_{\b}F^+_{\a^{-1}}-F^+_{\a}F^+_{\b^{-1}}
+F^+_{\a}F^+_{\b}F^+_{\a^{-1}}+F^+_{\b}F^+_{\a^{-1}}F^+_{\b^{-1}}
+F^+_{\a}F^+_{\b}F^+_{\a^{-1}}F^+_{\b^{-1}}.$$ \end{lemma} 
\proof We are going to use many times $F^+_{\gamma}=F_{\gamma}-1.$
We have
$$
\begin{tabular}{ll}
$F_{[\a,\b]}$&$=F_{\a}F_{\b}F_{\a^{-1}}F_{\b^{-1}}
=F_{\a}F^+_{\b}F_{\a^{-1}}F_{\b^{-1}}+F_{\a}F_{\a^{-1}}F_{\b^{-1}}=$\\
&$=F_{\a}F^+_{\b}F_{\a^{-1}}F_{\b^{-1}}+F_{\b^{-1}}
=F_{\a}F^+_{\b}F^+_{\a^{-1}}F_{\b^{-1}}+F_{\a}F^+_{\b}F_{\b^{-1}}+F_{\b^{-1}}=$\\
&$=F_{\a}F^+_{\b}F^+_{\a^{-1}}F_{\b^{-1}}+F_{\a}F_{\b}F_{\b^{-1}}-F_{\a}F_{\b^{-1}}+F_{\b^{-1}}=$\\
&$=F_{\a}F^+_{\b}F^+_{\a^{-1}}F_{\b^{-1}}+(F_{\a}-F_{\a}F_{\b^{-1}}+F_{\b^{-1}}-1)+1=$\\
&$=F_{\a}F^+_{\b}F^+_{\a^{-1}}F_{\b^{-1}}-F^+_{\a}F^+_{\b^{-1}}+1.$
\end{tabular}
$$
Also
$$
\begin{tabular}{ll}
$F_{\a}F^+_{\b}F^+_{\a^{-1}}F_{\b^{-1}}$&
$=F^+_{\a}F^+_{\b}F^+_{\a^{-1}}F_{\b^{-1}}+F^+_{\b}F^+_{\a^{-1}}F_{\b^{-1}}=$\\
&$=F^+_{\a}F^+_{\b}F^+_{\a^{-1}}F^+_{\b^{-1}}
+F^+_{\a}F^+_{\b}F^+_{\a^{-1}}+F^+_{\b}F^+_{\a^{-1}}F_{\b^{-1}}+F^+_{\b}F^+_{\a^{-1}}.$
\end{tabular}
$$
From these two sequences of equalities the lemma follows.


\subsection{Non-abelian reciprocity law}
Here we keep the notation from the previous section. Let $\omega_1\dots
\omega_n$ be differential forms of the third kind on a Riemann
surface $X$. Let $Y$ be the open subset of $X$ obtained by
removing the poles of $\omega_1\dots \omega_n$. Consider the
differential equation (\ref{DiffEq}) on $Y$
$$dF=F(\sum_{i=1}^n A_i\omega_i),$$
where $A_i$ for $i=1,\dots,n$ are non-commuting formal variables.
Fix a point $P$ in $Y$. Let
$$\gamma:[0,1]\rightarrow Y$$
be a piecewise smooth path that starts at $P$ and ends at $z$.
Then the solution of the differential equation with initial
condition $F(P)=1$ is
$$F_{\gamma}=1+\sum_{i=1}^n A_i\int_0^1\gamma^*\omega_i+
\sum_{i,j=1}^n A_iA_j\int_0^1\gamma^*\omega_i\circ\gamma^*\omega_j
+\dots.$$

If $\gamma_1$ and $\gamma_2$ are two paths such that the end point
of $\gamma_1$ is the beginning point of $\gamma_2$. Let,
$\gamma_1\gamma_2$ be the composition of the two paths. Then,
$F_{\gamma_1\gamma_2}=F_{\gamma_1}F_{\gamma_2}$ (theorem
\ref{thm1}).

Consider a simple loop $\sigma_{i}$ in $Y$ with the following
properties: it starts at $P$; it does not self-intersect; it
bounds a open region $V_{i}$ homeomorphic to a punctured disk; and
the only pole of the differential forms that lies inside the
closure of $V_i$ is $P_i$. Next we define a path $\gamma_i$ with
the following properties: $\gamma_i$ lies in the closure of $V_i$;
$\gamma_i$ starts at $P$ and ends at $P_i$. We can choose the
loops $\sigma_i$ so that in the counterclockwise we start with
loops around points of the poles of $\omega_1$. Then it continues
with loops around the points of the poles of $\omega_2$, which are
not poles of $\omega_1$, then it continues with loops around the
poles of $\omega_3$, which are not poles of $\omega_1$ or
$\omega_2$ and so on. Call these loops $\sigma_1,\dots,\sigma_N$.
We assume that $\sigma_i$ bounds a disk - the closure of $V_i$,
containing only one point $P_i$ from the poles of the differential
forms. Let also $\a_i$ and $\b_i$ be loops on $Y$ for
$i=1,\dots,g$, where $g$ is the genus of $Y$. We can choose such
that $\a_i$ and $\b_i$ so that $\a_i$ and $\b_i$ for $i=1,\dots,g$
and $\sigma_1,\dots,\sigma_N$ generate $\pi_1(Y)$ and the only
relation between them is
$$\sigma_1\dots\sigma_N[\a_1,\b_1]\dots[a_g,b_g]=1$$.

We define a pro-unipotent tame symbol to be $F_\sigma$ where
$\sigma$ is one of the loops that we have defined above. Note that
we have made some choices of loops. If we multiply all the loops
going in counterclockwise direction, then we will obtain a loop
homotopic to the zero loop at $P$. This gives the global
reciprocity law for pro-unipotent tame symbols after
simplification that uses section 1. We can write this reciprocity
law in the following way. Let $\sigma_1,\dots,\sigma_N$ in this
order be the above loops counted in counterclockwise direction.
Consider $F_{\sigma_i}$ as an element of the formal power series
${\mathbb C}<<A_1,\dots, A_n>>$.

In the formulation of the pro-unipotent reciprocity law, we are
going to use the notation from this section and from subsections
2.2 and 2.3.

\begin{theorem} \label{thm6} The non-abelian reciprocity law on a Riemann surface is
$$
\begin{tabular}{ll}
$(\prod_{i=1}^N
1+F^{reg}_{\gamma_i}F^{Res+}_{\sigma_i}F^{reg}_{\gamma_i^{-1}})\times$\\
$\times (\prod_{j=1}^g
1+F^+_{\b_j}F^+_{\a_j^{-1}}-F^+_{\a_j}F^+_{\b_j^{-1}}
+F^+_{\a_j}F^+_{\b_j}F^+_{\a_j^{-1}}+F^+_{\b_j}F^+_{\a_j^{-1}}F^+_{\b_j^{-1}}
+F^+_{\a_j}F^+_{\b_j}F^+_{\a_j^{-1}}F^+_{\b_j^{-1}})=1.$
\end{tabular}
$$ \end{theorem}
\proof For the solutions of this differential equation we have
that for two loops $\sigma$ and $\tau$ starting at the same point
$P$ we have that $F_{\sigma}F_{\tau}=F_{\sigma\tau}$ (see theorem
\ref{thm2}). The composition
$$\sigma_1\dots\sigma_N[\a_1,\b_1]\dots[\a_g,\b_g]$$
is homotopic to the trivial loop at $P$. Therefore,
$$F_{\sigma_1}\dots F_{\sigma_N}F_{[\a_1,\b_1]}\dots F_{\a_g,\b_g]}=1.$$
Using lemma \ref{lem5} in subsection 2.2, we have
$$F_{\sigma_i}=1+F^{reg}_{\gamma_i}F^{Res}_{\sigma_i}F^{reg}_{\gamma_i^{-1}}.$$
From theorem \ref{thm5} from subsection 2.3, we have
$$F_{[\a_j,\b_j]}=1+F^+_{\b_j}F^+_{\a_j^{-1}}-F^+_{\a_j}F^+_{\b_j^{-1}}
+F^+_{\a_j}F^+_{\b_j}F^+_{\a_j^{-1}}+F^+_{\b_j}F^+_{\a_j^{-1}}F^+_{\b_j^{-1}}
+F^+_{\a_j}F^+_{\b_j}F^+_{\a_j^{-1}}F^+_{\b_j^{-1}}.$$

The following two lemmas are useful for explicit computations,
which we are going to consider in section 4.
\begin{lemma} \label{lem7} The coefficient next to the linear terms in the formal
variables in
$$\prod_{i=1}^N F_{\sigma_{i}}$$ is zero.
\end{lemma}
\proof It is enough to prove the lemma for the coefficient
$A_1$. The coefficient next to $A_1$ is
$$\sum_{i=1}^N\int_{\sigma_{i}}\omega_1.$$
The sum is zero because it is equal to the sum of the residues of
$\omega_1$.

\begin{lemma} \label{lem8} The coefficient next to $A$, $B$ or $C$ in
$$F_{[\a_l,\b_l]}$$ is zero for every $l=1,\dots,g.$ \end{lemma}
\proof From lemma \ref{lem6} we see that there is no term in
degree 1 in the formal variables.

\section{Explicit formulas for reciprocity laws}
\subsection{Classical reciprocity laws}Let $\omega_1$ and $\omega_2$ be differential forms of
the third kind on a Riemann surface $X$ of genus at least $g\geq
1$. Consider the differential equation
$$dF=F(A\omega_1+B\omega_2),$$
(see the differential equation (\ref{DiffEq})) where $A$ and $B$
are non-commuting formal variables. Let $P$ be point on $X$, which
is not a pole for $\omega_1$  and $\omega_2$. Assume that there
are no common poles between $\omega_1$ and $\omega_2$.

Let $\gamma$ is a path starting at $P$ and ending at $Q$. Let
$F_\gamma=F(Q)$ be the solution of the differential equation with
initial conditions $F(P)=1$ solved along the path $\gamma$ and
$F(Q)$ is the evaluation of that solution at the point $Q$.

Let $P_1,\dots P_p$ be the poles of $\omega_1$. Let $Q_1,\dots
Q_q$ be the poles of $\omega_2$. Let
$Y=X-\{P_1,\dots,P_p,Q_1,\dots,Q_q\}$. Consider a simple loop
$\sigma_{P_i}$ in $Y$ with the following properties: it starts at
$P$; it does not self-intersect; it bounds a open region $V_{P_i}$
homeomorphic to a punctured disk; and the only pole of the
differential forms that lies inside the closure of $V_{P_i}$ is
$P_i$. We can choose these loops so that they do not intersect
each other except at the point $P$. Call these loops
$\sigma_{P_1},\dots,\sigma_{P_p},\sigma_{Q_1},\dots,\sigma_{Q_q}.$
We can choose the loops so that in the counterclockwise we start
with loops $\sigma_{P_1},\dots,\sigma_{P_p}$ around the poles of
$\omega_1$ in this order, followed by the loops
$\sigma_{Q_1},\dots,\sigma_{Q_q}$ around the poles of $\omega_2$.
Let also $\a_l$ and $\b_l$ be loops on $Y$ for $l=1,\dots,g$,
where $g$ is the genus of $Y$. We can choose $\a_l$ and $\b_l$ so
that $\a_l$ and $\b_l$ for $l=1,\dots,g$ and
$\sigma_{P_1},\dots,\sigma_{P_p},\sigma_{Q_1},\dots,\sigma_{Q_q}$
generate $\pi_1(Y,P)$ and the only relation between them is
$$\sigma_{P_1}\dots\sigma_{P_p}
\sigma_{Q_1}\dots\sigma_{Q_q}[\a_1,\b_1]\dots[\a_g,\b_g]=1.$$

Let $\gamma_{P_i}$ be a path from $P$ to $P_i$ that lies inside
the disk bound by $\sigma_{P_i}$. Similarly, let $\gamma_{Q_j}$ be
a path from $P$ to $Q_j$ that lies inside the disk bound by
$\sigma_{Q_j}$.

Let $\sigma$ be any of the above loops $\sigma_{P_i}$,
$\sigma_{Q_j}$, $\a_k$ or  $b_k$. Consider the tame symbol
$F_{\sigma}$.

Consider the coefficients contributing to $AB$ in the reciprocity
law associated to the differential equation
$$dF=F(A\omega_1+B\omega_2).$$

\begin{theorem} \label{thm7} Let $\omega_1$ and $\omega_2$ be two differential forms
of the third kind without common poles. Then with the above
notation we have a reciprocity law for differential forms of the
third kind
$$
\begin{tabular}{ll}
$\sum_{i=1}^pRes_{P_i}\omega_1\int_{\gamma_{P_i}^{-1}}\omega_2
+\sum_{j=1}^qRes_{Q_j}\omega_2\int_{\gamma_{Q_j}}\omega_1+$\\
$+\sum_{k=1}^g(\int_{\a_k}\omega_1\int_{\b_k}\omega_2
-\int_{\b_k}\omega_1\int_{\a_k}\omega_2)=$\\
$=0.$
\end{tabular}
$$
\end{theorem} \proof We have that
$$F_{\sigma_{P_1}\dots\sigma_{P_p}
\sigma_{Q_1}\dots\sigma_{Q_q}[\a_1,\b_1]\dots[\a_g,\b_g]}=1.$$
From lemma \ref{lem7}  and lemma \ref{lem8}  we have that there is
no linear terms in $A$ or $B$ in $F_{\sigma_{P_1}\dots\sigma_{P_p}
\sigma_{Q_1}\dots\sigma_{Q_q}}$ and in
$F_{[\a_1,\b_1]\dots[\a_g,\b_g]}$. Using theorem \ref{thm6}, lemma
\ref{lem6} and lemma \ref{lem5}, we obtain that the coefficient
next to $AB$ in $F_{\sigma_{P_i}}$ is
$$Res_{P_i}\omega_1\int_{\gamma_{P_i}^{-1}}\omega_2.$$ Similarly, the
coefficient next to $AB$ in $F_{\sigma_{Q_j}}$ is
$$Res_{Q_j}\omega_2\int_{\gamma_{Q_j}}\omega_1.$$ And finally, the
coefficient next to $AB$ in $F_{[\a_k,\b_k]}$ is
$$\int_{\a_k}\omega_1\int_{\b_k}\omega_2
-\int_{\b_k}\omega_1\int_{\a_k}\omega_2.$$


When the two forms do not have a common pole we obtain a
reciprocity law for differential forms of the third kind (see
\cite{GH}).

This also implies Weil reciprocity law. For similar approach see
\cite{Kh}.
 \begin{corollary} \label{co3}(Weil reciprocity) With the above notation,
consider $\omega_1=df/f$ and $\omega_2=dg/g$, where $f$ and $g$
are meromorphic functions on the Riemann surface with disjoint
divisors. Let the corresponding divisors $(f)$ and $(g)$ be
$$(f)=\sum_{i=1}^p a_i P_i$$
and
$$(g)=\sum_{j=1}^q b_j Q_j.$$
Then
$$
\prod_{i=1}^p g(P_i)^{-a_i} \prod_{j=1}^q f(Q_j)^{b_j}=1.
$$
\end{corollary} \proof Consider the above theorem with $\omega_1=df/f$ and
$\omega_2=dg/g$ and exponentiate.


\subsection{Non-abelian reciprocity laws involving three
differential forms}
Consider the differential equation (\ref{DiffEq}) in the following
case
$$dF=F(A\omega_1+B\omega_2+C\omega_3),$$
where $A$, $B$ and $C$ are non-commuting formal variables. Assume
$\omega_1$, $\omega_2$ and $\omega_3$ have distinct poles. Let $P$
be point on $X$, which is not a pole for $\omega_1$, $\omega_2$ or
$\omega_3$.

Let $\gamma$ is a path starting at $P$ and ending at $Q$. Let
$F_\gamma=F(Q)$ be the solution of the differential equation with
initial conditions $F(P)=1$ solved along the path $\gamma$ and
$F(Q)$ is the evaluation of that solution at the point $Q$.

Let $P_1,\dots P_p$ be the poles of $\omega_1$, $Q_1,\dots Q_q$ be
the poles of $\omega_2$ and $R_1,\dots R_r$ be the poles of
$\omega_3$. Let
$Y=X-\{P_1,\dots,P_p,Q_1,\dots,Q_q,R_1,\dots,R_r\}$. Consider a
simple loop $\sigma_{P_i}$ in $Y$ with the following properties:
it starts at $P$; it does not self-intersect; it bounds a open
region $V_{P_i}$ homeomorphic to a punctured disk; and the only
pole of the differential forms that lies inside the closure of
$V_{P_i}$ is $P_i$. We can choose these loops so that they do not
intersect each other except at the point $P$. Call these loops
$\sigma_{P_1},\dots,\sigma_{P_p},\sigma_{Q_1},\dots,\sigma_{Q_q},
\sigma_{R_1},\dots,\sigma_{R_r}$.
We can choose the loops so that in the counterclockwise we start
with loops $\sigma_{P_1},\dots,\sigma_{P_p}$ around the poles of
$\omega_1$ in this order, followed by the loops
$\sigma_{Q_1},\dots,\sigma_{Q_q}$ around the poles of $\omega_2$,
followed by the loops $\sigma_{R_1},\dots,\sigma_{R_r}$ around the
poles of $\omega_3$. Let also $\a_l$ and $\b_l$ be loops on $Y$
for $l=1,\dots,g$, where $g$ is the genus of $Y$. We can choose
$\a_l$ and $\b_l$ so that $\a_l$ and $\b_l$ for $l=1,\dots,g$ and
$\sigma_{P_1},\dots,\sigma_{P_p},\sigma_{Q_1},\dots,
\sigma_{Q_q},\sigma_{R_1},\dots,\sigma_{R_r}$
generate $\pi_1(Y,P)$ and the only relation between them is
$$\sigma_{P_1}\dots\sigma_{P_p}\sigma_{Q_1}\dots\sigma_{Q_q}
\sigma_{R_1}\dots\sigma_{R_r}[\a_1,\b_1]\dots[\a_g,\b_g]=1$$.

Let $\gamma_{P_i}$ be a path from $P$ to $P_i$ that lies inside
 the punctured disk $V_{P_i}$ bound by $\sigma_{P_i}$. Similarly, let
$\gamma_{Q_j}$ be a path from $P$ to $Q_j$ that lies inside the
punctured disk $V_{Q_j}$ bound by $\sigma_{Q_j}$ and let
$\gamma_{R_k}$ be a path from $P$ to $R_k$ that lies inside the
punctured disk $V_{R_k}$ bound by $\sigma_{R_k}$

Let $\sigma$ be any of the above loops $\sigma_{P_i}$,
$\sigma_{Q_j}$, $\sigma_{R_k}$ $\a_l$ or  $b_l$. Consider the tame
symbol $F_{\sigma}$. For them we have the following reciprocity.

\begin{theorem} \label{thm8} Using the above notation, for three differential forms
of the third kind, with disjoint poles, we have
$$
\begin{tabular}{l}
$\sum_{i=1}^p
Res_{P_i}\omega_1\int_{\gamma_{P_i}^{-1}}\omega_2\circ\omega_3+
\sum_{j=1}^q
Res_{Q_j}\omega_2\int_{\gamma_{Q_j}^{-1}}\omega_1\int_{\gamma_{Q_j}}\omega_3+$\\
$\sum_{k=1}^r
Res_{R_k}\omega_3\int_{\gamma_{R_k}}\omega_1\circ\omega_2+$\\
$+\sum_{l=1}^g
(\int_{\a_l}\omega_1\circ\omega_2\int_{\b_l}\omega_3-
\int_{\b_l}\omega_1\circ\omega_2\int_{\a_l}\omega_3+$\\
$+\int_{\a_l}\omega_3\circ\omega_2\int_{\b_l}\omega_1-
\int_{\b_l}\omega_3\circ\omega_2\int_{\a_l}\omega_1-$\\
$-\int_{\a_l}\omega_1\int_{\b_l}\omega_2\int_{\a_l}\omega_3+
\int_{\b_l}\omega_1\int_{\a_l}\omega_2\int_{\b_l}\omega_3)=0.$
\end{tabular}
$$
\end{theorem} \proof

Consider only the coefficients of the tame symbol $F_{\sigma}$
next to the $A$, $B$ $C$, $AB$, $BC$ and $ABC$ terms in the
non-commuting power series in $A$ $B$ and $C$.

The coefficient next to $ABC$ in $F_{\sigma_{P_i}}$ is
$$\int_{\sigma_{P_i}}\omega_1\circ\omega_2\circ\omega_3=
Res_{P_i}\omega_1\int_{\gamma_{P_i^{-1}}}\omega_2\circ\omega_3.$$

Let $\sigma_1=\sigma_{P_1}\sigma_{P_2}\dots\sigma_{P_p}$ be a
product of loops that go around the poles of $\omega_1$. Then the
coefficient next to $ABC$ in $F_{\sigma_1}=F_{\sigma_{P_1}}\dots
F_{\sigma_{P_p}}$ is
$$\int_{\sigma_1}\omega_1\circ\omega_2\circ\omega_3=
\sum_{i=1}^p
Res_{P_i}\omega_1\int_{\gamma_{P_i}^{-1}}\omega_2\circ\omega_3.$$
Also the coefficient next to $A$ in $F_{\sigma_1}$ is equal to the
sum of the residues of $\omega_1$, which is zero. Also the
coefficient next to $B$ and $C$ in $F_{\sigma_1}$ is zero because
the disks that $\sigma_1$ bounds do not contain any of the poles
of $\omega_2$ or $\omega_3$.

Let $\sigma_3=\sigma_{R_1}\sigma_{R_2}\dots\sigma_{R_r}$ be the
product of the loops around the poles of $\omega_3$. Consider the
coefficients next to $A$, $B$ and $C$ in $F_{\sigma_3}$. The
coefficient next to $A$ and $B$ is zero because the loop does not
bound any poles of $\omega_1$ or $\omega_2$. The coefficient next
to $B$ is zero because the sum of the residues of $\omega_3$ is
zero. Similarly to the commutation in the previous paragraph the
coefficient next to $ABC$ in $F_{\sigma_3}$ is
$$\int_{\sigma_3}\omega_1\circ\omega_2\circ\omega_3=
\sum_{k=1}^r
Res_{R_k}\omega_3\int_{\gamma_{R_k}}\omega_1\circ\omega_2.$$

Let $\sigma_2=\sigma_{Q_1}\sigma_{Q_2}\dots\sigma_{Q_q}$ be the
product of the loops around the poles of $\omega_2$. Consider the
coefficients next to $A$, $B$ and $C$ in $F_{\sigma_2}$. The
coefficient next to $A$ and $C$ is zero because the loop does not
bound any poles of $\omega_1$ or $\omega_3$. The coefficient next
to $B$ is zero because the sum of the residues of $\omega_2$ is
zero.

The coefficient next to $ABC$ in $F_{\sigma_2}$ is
$$\int_{\sigma_2}\omega_1\circ\omega_2\circ\omega_3=
\sum_{j=1}^q
Res_{Q_j}\omega_2\int_{\gamma_{Q_j}^{-1}}\omega_1\int_{\gamma_{Q_j}}\omega_3.$$

Therefore, the coefficient next to $ABC$ in
$F_{\sigma_1\sigma_2\sigma_3}$ is
$$\sum_{i=1}^p
Res_{P_i}\omega_1\int_{\gamma_{P_i}^{-1}}\omega_2\circ\omega_3+
\sum_{j=1}^q
Res_{Q_j}\omega_2\int_{\gamma_{Q_j}^{-1}}\omega_1\int_{\gamma_{Q_j}}\omega_3+
\sum_{k=1}^r
Res_{R_k}\omega_3\int_{\gamma_{R_k}}\omega_1\circ\omega_2.$$

Consider the coefficients next $A$, $B$, $C$ and $ABC$ of
$F_{[\a_l,\b_l]}.$ Using lemma \ref{lem6}, we notice that the
coefficient next to $A$, $B$ or $C$ in $F_{[\a_l,\b_l]}$ is zero.
Also, the linear coefficients next to $F_{\sigma_1}$,
$F_{\sigma_2}$ and $F_{\sigma_3}$ are zero. Therefore, there will
be any contribution from the quadratic terms in $A$, $B$ and $C$
in the series $F_{[\a_l,\b_l]}$, since it has to be multiplies by
a linear term, which is zero. From the reciprocity law, theorem
\ref{thm6}, we have consider only the coefficient next to $ABC$
of. We only need to consider the sum of the coefficients next to
$ABC$ in $F_{[\a_l,\b_l]}$, $F_{\sigma_1}$, $F_{\sigma_2}$ and
$F_{\sigma_3}$, because the linear terms are zero. By lemma
\ref{lem8} (see also theorem \ref{thm6}), we have that the
coefficient next to $ABC$ in $F_{[\a_l,\b_l]}$ is
$$
\begin{tabular}{ll}
$\int_{\b_l}\omega_1\int_{\a_l^{-1}}\omega_2\circ\omega_3 +
\int_{\b_l}\omega_1\circ\omega_2\int_{\a_l^{-1}}\omega_3-$\\
$-\int_{\a_l}\omega_1\int_{\b_l^{-1}}\omega_2\circ\omega_3 -
\int_{\a_l}\omega_1\circ\omega_2\int_{\b_l^{-1}}\omega_3+$\\
$+\int_{\a_l}\omega_1\int_{\b_l^{-1}}\omega_2\int_{\a_l^{-1}}\omega_3+
\int_{\b_l}\omega_1\int_{\a_l^{-1}}\omega_2\int_{\b_l^{-1}}\omega_3$
\end{tabular}
$$
Using lemma \ref{lem1} for reversing a path of the above iterated
integrals, we finish the proof of the theorem.
\begin{remark}
We need the last formula, for the coefficient next to $ABC$ in
$F_{[\a_l,\b_l]}$, in the paper on a refinement of the Parshin
symbol for surfaces \cite{H}.
\end{remark}
\section{Manin's non-commutative modular symbol, involving Eisenstein series}
In this section we give an application of our reciprocity law to
the non-commutative modular symbol that Manin has defined (see
\cite{M}). Before we can apply the reciprocity law, we make a
generalization of Manin's non-commutative modular symbol. We
generalize it so that we can consider modular symbols involving
Eisenstein series.

\subsection{Construction and properties of the symbols}
In this section we enrich Manin's non-commutative modular symbol
so that it includes Eisenstein series, not only cusp forms. This
was a question raised by Glen Stevens.

Let $\Gamma$ be an arithmetically defined, torsion free, discrete
group, acting on the upper half plane ${\cal{H}}$, whose quotient
is  a Riemann surface $X$ without a divisor $D$. (We allow compact
quotients. That is $D=0$.) We will consider loops on the modular
curve $X-D$. Let $f_1dz,\dots,f_ndz$ be cusp or Eisenstein modular
forms on $\cal{H}$ with respect to $\Gamma$. Consider the
differential equation
\begin{equation} dJ(z)=J(z)\sum_{i=1}^n A_i f_idz,
\label{ModSymb}\end{equation} where $A_i$'s are non-commuting
formal variables. This is an equation of the type of equation
(\ref{DiffEq}), which Manin uses to define the non-commutative
modular symbol \cite{M}. Originally Manin considered only the case
where all $f_i$ are cup forms. Here we allow them to be either
cusp forms or Eisenstein series.

Let $P$ be the image of $i=\sqrt{-1}$ in the projection
$$pr:{\cal{H}}\rightarrow X-D.$$
That is, $$P=pr(\sqrt{-1}) .$$

Usually, by a modular form $f(z)$ people denote a holomorphic {\it
function} on the upper half plane that has 'good' transformation
properties under the action of the arithmetic group $\Gamma$. W
are going to use a different trivialization of the modular form. A
modular form will be a holomorphic {\it $1$-form} on the upper
half plane that has a 'good' transformation properties under the
action of the group $\Gamma$.

We are going to identify the upper half plane ${\cal H}$ with
$SL_2({\mathbb R})/SO_2({\mathbb R})$ by identifying $i$ in the
upper half plane ${\cal H}$ with the projection of the identity
element $I$ in $SL_2({\mathbb R})$. Let $P\in X-D$ be the
projection of $i$ from the upper half plane. Then
$$\pi_1(X-D,P)=\Gamma.$$
Let $Q_1,\dots,Q_N$ be the points of the divisor $D$. They are
called cusp points. Consider the punctured Riemann surface $X$. We
can slice it to $(N-1+4g)$-gon, so that the points of the divisor
$D$ are vertices of this polygon.

Denote by $\gamma_j$ a path, which is a piece of one of the sides
of this polygon, that starts at an orbit of $P$ and ends at a
vertex with label $Q_j$. Denote by $\tilde{\gamma}_j^{-1}$ a path,
which is a piece of one of the sides of this polygon, that starts
at $Q_j$ and ends at the other orbit of $P$ compared to
$\gamma_j$.

Let $\omega_1,\dots,\omega_n$ be modular $1$-forms. Consider a
neighborhood $U$ of a cusp point $Q=Q_j$ on the modular curve $X$
such that the modular $1$-forms $\omega_1,\dots,\omega_n$ have
trivialization on $U$. Let $Q'$ be in $U$ and $Q'\neq Q$. As
before, let $\sigma=\sigma_j$ be a simple loop on the modular
curve $X-D$ that starts at $P$ and bounds a disc, which contains
$Q$. Let $Q_\epsilon$ be a point inside $U$ in an
$\epsilon$-neighborhood of $Q$. The loop $\sigma$ is homotopic to
$$\gamma'\gamma_\epsilon\sigma_e\gamma_e^{-1}\gamma'^{-1},$$
where $\gamma'$ is a path from $P$ to $Q'$, $\gamma_e$ a path from
$Q'$ to $Q_\epsilon$ and$\sigma_e$ a loop, starting at $Q_e$
around $Q$ in an $\epsilon$-neighborhood of $Q$. Let
$\tilde{\sigma}=\tilde{si}_j$ be a lift of $\sigma$ on the upper
half plane. Define lifts $\tilde{\gamma}'$,
$\tilde{\gamma}_\epsilon$, $\tilde{\sigma}_e$,
$\tilde{\gamma}_\epsilon^{(-1)}$, $\tilde{\gamma}'^{(-1)}$ from
the modular curve $X-D$ to the upper half plane of $\gamma'$,
$\gamma_\epsilon$, $\sigma_e$, $\gamma_e^{-1}$ and $\gamma'^{-1},$
respectively. We also require that the ending point of a lifted
path is the staring point of the next path.

\begin{definition}
\label{def3} Given modular $1$-forms
$(\omega_1,\dots,\omega_n)=\Omega$ and a path $\gamma$ in the
upper half plane ${\cal H}$, define
$$J_\gamma(\Omega)$$ to be the generating series of iterated
integrals of $\omega_1,\dots,\omega_n$ over the path $\gamma$, if
the modular forms do not have poles at the starting point and at
the ending point. If the starting point of $\gamma$ is a pole of
any of the modular $1$-forms, then we define
$$J^{Reg}_\gamma(\Omega)$$ to be the generating series of iterated
integrals whose first $1$-form in the iteration has no pole at the
starting point of $\gamma$. Similarly, if the ending point of
$\gamma$ is a pole of any of the modular $1$-forms, then we define
$$J^{Reg}_\gamma(\Omega)$$ to be the generating series of iterated
integrals whose last $1$-form in the iteration has no pole at the
ending point of $\gamma$. We define
$$J^{Res}_Q(\Omega)$$
where $Q$ is a cusp point, to be exponential of the sum of the
residues of the modular $1$-forms at the point $Q$ on the modular
curve, computed in a small neighborhood of $Q$. And finally,
$$J^{Res+}_Q(\Omega)=-1+J^{Res}_Q(\Omega).$$ \end{definition}

We have
\begin{theorem} \label{thm11} With the notation from this subsection, we have
$$J_{\tilde{\sigma}}(\Omega)=1+J^{reg}_{\tilde{\gamma}}
(\Omega)J^{Res+}_{Q}(\Omega)J^{reg}_{\tilde{\gamma}^{(-1)}}(\Omega),$$
where $\gamma$ is a path starting at $P$  and ending at $Q$,
$\tilde{\gamma}$ and $\tilde{\gamma}^{(-1)}$ are lifts of $\gamma$
and $\gamma^{-1}$, so that $\tilde{\gamma}$ starts at the initial
point of $\tilde{\sigma}$ and ends at a cusp corresponding to $Q$
and $\tilde{\gamma}^{(-1)}$ at the ending point of
$\tilde{\gamma}$ and ends at the ending point of
$\tilde{\sigma}$\end{theorem}
\proof We apply theorem \ref{thm5} to the loop
$\gamma_e\sigma_e\gamma_e^{-1}$.Then we use that
$J_{\gamma'}(\Omega)J_{\gamma_e}(\Omega)=J_{\gamma'\gamma_e}(\Omega)$.
The same equality hold for regularized generating series.

We have a non-abelian reciprocity law on the upper half plane
involving holomorphic modular $1$-forms, which are either cusp
forms or Eisenstein series (as $1$-forms). It is similar to the
non-abelian reciprocity law on Riemann surface involving
differential forms of he third kind. See theorem \ref{thm6}. In
the following theorem we are going to omit $\Omega$ from the
notation. That is, we are going to write $J_{\tilde{\gamma}}$
instead of $J_{\tilde{\gamma}}(\Omega)$.
\begin{theorem} \label{thm12}With the above notation, we have the following
reciprocity law for Eisenstein and cusp forms
$$\left(\prod_{j=1}^M
1+J^{reg}_{\tilde{\gamma}_j}J^{Res+}_{Q_j}J^{reg}_{\tilde{\gamma}_j^{(-1)}}\right)=1,
$$
where the composition of the paths $\tilde{\gamma}_j$ and
$\tilde{\gamma}_j^{(-1)}$ for all $j$, namely,
$$\prod_{j=1}^M\tilde{\gamma}_j\tilde{\gamma}_j^{(-1)}$$
is contractible and $M=N+4g$. \end{theorem}
\begin{definition}
\label{def4} Let $P$ and $P'$ be in the orbit of $i$ in the upper
half plane under the action of the arithmetic group $\Gamma$. Let
$Q$ and $Q'$ be a cusp points. Let also $\tilde{\gamma}$,
$\tilde{\gamma}'^{(-1)}$ and $\tilde{\a}$ be paths in the upper
half plane that connect $P$ with $Q$, $Q'$ with $P$ and $P$ with
$P'$, respectively. We define  non-commutative modular symbols to
be
$$J^{reg}_{\tilde{\gamma}}(\Omega),\mbox{ }J^{reg}_{\tilde{\gamma}'^{(-1)}}(\Omega),
\mbox{ }J_{\tilde{\a}}(\Omega)$$\mbox{ and }
$$J^{Res}_{Q}(\Omega).$$
\end{definition} The non-commutative symbol that Manin has defined can be
recovered by considering the product
$$J^{reg}_{\tilde{\gamma}'^{(-1)}}(\Omega)J^{reg}_{\tilde{\gamma}}(\Omega),$$
where the regularization is not needed for cusp forms.
\subsection{Main Example}
A very important example occurs, when we consider generating
series of iterated integrals $J$ of $1$-forms $fdz$ and $dz$ for a
torsion free arithmetic groups $\Gamma$ commensurable to
$SL_2({\mathbb Z})$, where $fdz$ is a cusp form.

\begin{remark} Note that $fdz$ is a cusp $1$-form but $f$ is not a cusp
form. The poles of $f$ cancel with the poles of $dz$. In terms of
algebraic geometry a cusp form is a section of a certain line
bundle on the modular curve $X$. Instead of considering local
trivialization in terms of the rational functions ${\cal O}_X$, we
consider local trivialization in terms of $1$-forms on $X$,
$\Omega_X$. The cusp form $F$ trivialized as a function on the
upper half plane is related to the cusp form trivialized as a
holomorphic $1$-form $fdz$  by
$$F(x)=\int_0^x fdz$$
for a point $x$ in the chosen fundamental domain.
\end{remark}
\begin{definition}
\label{def5} We denote by $$dz^{\circ n}$$ the $n$-fold iteration
of the form $dz$. \end{definition}
With this notation we have the following theorem.
\begin{theorem} \label{thm13}
Let $$L(F,n_1,F,n_2,\dots,F,n_k):=\int_{i\infty}^0 F(z)
z^{n_1-1}dz\circ\dots\circ F(z)z^{n_k-1}dz.$$ Then

$$\begin{tabular}{ll}(a)& $L(F,n)=(n+1)!\int_{i\infty}^0 fdz\circ dz^{\circ n};$\\
\\
(b)&$L(F,n_1,F,n_2,\dots,F,n_k)=\prod_{j=1}^k(n_j+1)!\int_{i\infty}^0
fdz\circ dz^{\circ n_1}\circ\dots\circ fdz\circ dz^{\circ n_k};$\\
\\
(c)&$J^{reg}_{\tilde{\gamma}}(fdz,dz)$ is the generating series of
the $L$-functions in part (b).
\end{tabular}$$
\end{theorem}
\proof The form $dz$ on the upper half plane corresponds to a
differential form of the third kind on the modular curve. Also, at
the point $0$ of the upper half plane, both $fdz$ and $dz$ are
holomorphic. Thus, the regularization is well defined.

{\bf Acknowledgments.} I would like to thank Yuri Manin for his
inspiring talk on non-commutative modular symbol. Also, I would
like to thank Alexander Goncharov, from whom I learned iterated
integrals. And finally, I would like to thank Glen Stevens for
raising the question about non-commutative modular symbols
involving Eisenstein series.

I would like to thank Max Planck Insitut f\"ur Mathematik in Bonn,
the University of Durham, Brandeis University and Universit\"at
T\"ubinbgen for the kind hospitality and the creative atmosphere.
Also, I would like to thank the Arithmetic Algebraic Geometry
Marie Curie Network for the financial support.





\end{document}